\documentclass[11pt]{article}   % For Latex2e

\usepackage{amsthm, amsmath, amssymb,epsfig}   % For Latex2e

\newcommand{\oa}{\overline{a}}
\newcommand{\ox}{\overline{x}}
\newcommand{\olam}{\overline{\lambda}}

\newcommand{\lpp}{\operatorname{lpp}}

\newcommand{\lc}{\operatorname{lc}}

\newcommand{\lex}{\operatorname{lex}}
\newcommand{\cont}{\operatorname{cont}}

\newcommand{\V}{\mathbb V}
\newcommand{\I}{\mathbb I}
\newcommand{\hs}{\hspace*{3mm}}
\newcommand{\Compl}{\mathbb C}
\newcommand{\Rac}{\mathbb Q}
\newcommand{\Ent}{\mathbb Z}
\newcommand{\anul}[1]{}

\newtheorem{theorem}{Theorem}%[section]
\newtheorem{lemma}[theorem]{Lemma}

\newtheorem{conjecture}[theorem]{Conjecture}

\theoremstyle{definition}
\newtheorem{definition}[theorem]{Definition}

\title{On the canonical discussion of polynomial systems with parameters\footnote{Work partially supported by the
    Ministerio de Ciencia y Tecnolog\'ia under project
    MTM 2006-01267, and by the Generalitat de Catalunya under project
    2005 SGR 00692}}

\author{Antonio Montes\\ \\
Departament de Matem\`atica Aplicada 2,\\ Universitat
Polit\`ecnica de Catalunya, Spain.\\
e-mail: antonio.montes@upc.edu\\
http://www-ma2.upc.edu/$\sim$montes}

\date{Mai 2007}

\begin{document}
\anul{
\titlepage{
\null \vspace{2.5truecm} \thispagestyle{empty}
\thispagestyle{empty}
\begin{center}
\textbf{\large On the canonical discussion of \\
polynomial systems with parameters}
%\textbf{\large Ascending discriminant chains for the \\ discussion of parametric polynomial ideals}
\\[3mm]
\textsc{Antonio Montes}
 \vspace{0.8truecm}

MA2--IR--06--00006
\end{center}
}
\newpage
\null \thispagestyle{empty}
\newpage

\setcounter{page}{1} }

 \maketitle
%Begin new own definitions

\begin{abstract}
  Given a parametric polynomial ideal $I$, the algorithm DISPGB,
  introduced by the author in 2002, builds up a binary tree
  describing a dichotomic discussion of the different reduced
  Gr\"obner bases depending on the values of the parameters, whose set
  of terminal vertices form a Comprehensive Gr\"obner System (CGS).
  It is relevant to obtain CGS's having further
  properties in order to make them more useful for the applications.
  In this paper the interest is focused on obtaining a canonical CGS.
  We define the objective, show the difficulties and formulate a natural
  conjecture. If the conjecture is true
  then such a canonical CGS will exist and can be computed. We also give an
  algorithm to transform our original CGS in this direction and
  show its utility in applications.

\begin{tabular}{p{13.5cm}}
\end{tabular}

\end{abstract}

\section{Introduction}

There are many
authors~\cite{Be94,BeWe93,De99,DoSeSt06,Du95,FoGiTr01,Gi87,Gom02,
GoTrZa00,GoTrZa05,HeMcKa97,Ka97,Kap95,MaMo06,Mo02,Mor97,Pe94,
SaSu03,SuSa06,Sc91,Si92,We92,We03,Wi06} who have studied the
problem of specializing parametric ideals into a field and
determining the specialized Gr\"obner bases. Many other
authors~\cite{Co04,Em99,GoRe93,GuOr04,Mo95,Mo98,Ry00} have applied
some of these methods to solve concrete problems. In the previous
paper \cite{Mo02} we give more details of their contributions to
the field. In the following we only refer to the papers directly
related to the present work.

Let $I \subset K[\oa][\ox]$ be a parametric ideal in the variables
$\ox=x_1,\dots,x_n$ and the parameters $\oa=a_1,\dots,a_m$, and
$\succ_{\ox}$ and  $\succ_{\oa}$ monomial orders in variables and
parameters respectively. Denote $A=K[\oa]$.
Weispfenning~\cite{We92} proved the existence of a {\em
Comprehensive Gr\"obner Basis} (CGB) of $I$ and gave an algorithm
for computing it. It exists a modern implementation in REDUCE of
CGB algorithm due to T. Sturm et al.~\cite{DoSeSt06}. Let $K$ be a
computable field (for example $\Rac$) and ${K'}$ an algebraically
closed extension (for example $\Compl$). A CGB of $I\subset A[{\bf
x}]$ wrt (with respect to) the termorder $\succ_{{\bf x}}$ is a
basis of $I$ that specializes to a Gr\"obner basis of
$\sigma_{\oa_0}(I)$ for any specialization
$\sigma_{\oa_0}:K[\oa][\ox]\rightarrow {{K'}}[\ox]$, that
substitutes the parameters by values $\oa_0 \in {K'}^m$.

In most applications of parametric ideals the related object
called {\em Comprehensive Gr\"obner System} (CGS) is more
suitable. A CGS of the ideal $I\subset A[{\bf x}]$ wrt
$\succ_{{\ox}}$ is a set
\[
\begin{array}{lcl}
\operatorname{CGS}(I,\succ_{{\ox}})&=& \{ (S_i,B_i)\ :\ 1 \le i
\le s,\ S_i \subset K'^m,\ B_i
\subset A[{\ox}],\ \bigcup_i S_i=K'^m,  \\
&&  \forall {\oa_0} \in S_i,\  \sigma_{{\oa_0}}(B_i) \hbox{ is a
Gr\"obner basis of } \sigma_{{\oa_0}}(I) \hbox{ wrt }
\succ_{{\ox}} \}.
\end{array}
\]
The sets $S_i$ are often called {\em segments} and it is always
assumed that they are {\em constructible sets}. A CGB is a special
CGS with a unique segment $K'^m$. In a CGB the polynomials in the
basis are {\em faithful}, i.e. they belong to $I$. Further
properties are required to obtain more powerful CGS.
\begin{definition}[Disjoint CGS]\label{DefDisCGS}
 A CGS is said to be {\em disjoint} if the sets $S_i$ form a partition of
 $K'^m$.
 \end{definition}
\begin{definition}[Reduced basis]\label{DefRedBas}
A subset $B\subset A[\ox]$ is a reduced basis for a segment $S$ if
it verifies the following properties:\\
\begin{tabular}{rp{5in}}
  (i) & the polynomials in $B$ are normalized to have
  content 1 wrt $\ox$ over $A$ {\rm (}in order to work with polynomials instead of
  rational functions{\rm )};\\
  (ii) & the leading coefficients of the polynomials in $B$
  are different from zero on every point of $S$;\\
  (iii)& $B$ specializes to the
  reduced Gr\"obner basis of $\sigma_{\oa_0}(I)$, keeping the same
 $\lpp$ {\rm (}leading power product set{\rm )} for each ${\oa_0} \in S$,
 i.e. its $\lpp$ set remains stable under specializations within
 $S$.\\
\end{tabular}
\end{definition}
Reduced bases are {\em not faithful}, i.e. they do not, in
general, belong to $I$. They are not unique for a given segment,
but the number of polynomials as well as the $\lpp$ are unique.
\begin{definition}[Reduced CGS]\label{DefRedCGS}
 A CGS is said to be {\em reduced} if its segments have reduced
 bases.
\end{definition}
As it is known, the $\lpp$ of the reduced Gr\"obner basis of an
ideal determine the cardinal or dimension of the solution set over
an algebraically closed field. This is the reason why disjoint
reduced CGS are very useful for applications as they characterize
the different kind of solutions of $\V(I)$.

Using Weispfenning's suggestions the author~\cite{Mo02} obtained
an efficient algorithm (DISPGB) for Discussing Parametric
Gr\"{o}bner Bases to compute a disjoint reduced CGS. Actually this
algorithm is called BUILDTREE.

BUILDTREE builds up a dichotomic binary tree, whose branches at
each vertex correspond to the annihilation or not of a polynomial
in $K[\oa]$. It places at each vertex $v$ a specification
$\Sigma_v=(N_v,W_v)$ of the included specializations, that
summarizes the null and non-null decisions taken before reaching
$v$, and a specialized basis $B_v$ of $\sigma_{\oa}(I)$ for the
specializations $\sigma_{\oa_0} \in \Sigma_v$. The set of terminal
vertices form a disjoint reduced CGS where the segments $S_v$ are
characterized by reduced specifications determined by $(N_v,W_v)$.

Since then, more advances have been made. Inspired by BUILDTREE,
Weispfenning~\cite{We03} gave a constructive method for obtaining
a canonical CGB (CCGB) for parametric polynomial ideals.
%His algorithm is based on the direct
%computation of a discriminant ideal. Nevertheless his method does
%not provide a true partition and thus we cannot obtain from it a
%disjoint CGS.
Using this idea, Manubens and Montes~\cite{MaMo06} improved
BUILDTREE showing that the tree $T_0$ built up by BUILDTREE, can
be rewritten as a new tree $T$ providing a more compact and
effective discussion by computing a discriminant ideal that is
easy to compute from $T_{0}$. The rebuilding algorithm given in
~\cite{MaMo06} can be iterated to obtain a very compact new tree
organized as a right-comb tree. It builds an ascending chain of
discriminant ideals that orders the segments defined as
differences of the varieties of two consecutive discriminants, and
provides a very compact disjoint reduced CGS. Nevertheless this
rebuilding algorithm does not always produce the canonical CGS.
This method will be presented in a forthcoming paper.

Having in mind the improvement of our disjoint reduced CGS to
obtain a canonical CGS, in the present paper we adopt a different
perspective. Instead of rebuilding the tree, we analyze all the
different situations that can occur in the BUILDTREE CGS,
formulate a natural conjecture and, using it, show how the
segments can be packed to obtain the largest possible segments
allowing the same reduced basis. These segments become
non-algorithmic dependent, i.e. intrinsic for the given ideal. The
packed intrinsic partition contains the minimum number of segments
corresponding to reduced bases. Algorithms to perform the
discussion about which segments must be packed, to obtain the
reduced basis for the packed segments and to describe segments
defined by difference of two varieties in a canonical form are
also given. It remains to give a canonical description of the
union of the segments included in the packed segments as well as
the algorithm to carry it out. This last step is described
in~\cite{MaMo07a}.

The whole set of algorithms, including those in~\cite{MaMo07a}
have been implemented\footnote{Manubens and Montes implementation
of the new algorithm MCCGS is available on the web
http://www-ma2.upc.edu/$\sim$montes. The library, called DPGB
release 7, is implemented in {\em Maple} 8.} and denoted MCCGS
(Minimal Canonical Comprehensive Gr\"obner System) algorithm. It
is yet operative and in experimental phase. It is promising and
very useful for applications as can be seen through its
applications to automatic geometric theorem proving and
discovering~\cite{MoRe07}. It may be objected that canonicity pays
a price in computing time. The implementation shows that in fact
the time increases only about 20-30\% whether the output becomes
much more simpler, compact, easy to be understood and practical
for applications.

Recently Sato and Suzuki~\cite{SuSa06} have developed a new simple
algorithm for computing CGS based on Kalkbrenner's
Theorem~\cite{Ka97}. Its interest lies in its simplicity (it is
perhaps sometimes more efficient) but the output is not
sufficiently clear and useful for applications.

Section \ref{SecBuildtree} reviews the basic features of
BUILDTREE. In Section \ref{SecFindCan} it is explained what is
meant by canonical CGS, a conjecture is formulated and using it,
it is shown how to obtain the canonical CGS.
Section~\ref{SecCanSpec} give the basic theorems to obtain a
canonical description of diff-specifications and gives the
corresponding algorithms. Finally, in Section~\ref {SecFurther}
further developments are pointed out giving some insight as an
advance of the content of the final paper~\cite{MaMo07a}, where
the whole canonical description of the union of segments will be
given making the MCCGS the algorithm to provide a canonical
representation of the intrinsic segments.

In this paper, we only give partial examples that illustrate the
algorithms discussed here. A unique complete example is given in
the final Section~\ref{SecFurther} to give an idea of how simple
is the final output that has the minimum number of segments.

The algorithms work with ideals, as these are the algebraic
objects that allow a Gr\"obner representation. But ideals do not
represent varieties in a unique form. So we frequently adopt a
geometrical view. We shall consider ideals defined in $A=K[\oa]$,
where $K$ is the computable field (for example $\Rac$), whereas
the varieties will be considered in ${K'}^m$, where $K'$ is an
algebraically closed extension (for example $\Compl$). Let $J$ be
an ideal in $K[\oa]$, and $J'=J\cdot K'[\oa]$ be its extension to
$K'[\oa]$. The symbol $\V(J)$ will denote
\[
\begin{array}{lcl}
\V(J)&=&\{\oa \in {K'}^m : \forall f \in J, f(\oa)=0 \} \\
&=& \{\oa \in {K'}^m : \forall f \in J', f(\oa)=0 \}=\V(J').
\end{array}
\]
We emphasize the use of the non-standard notation $\V$ in the
whole paper as used in the extended affine space, whereas the
ideals are defined in the base field. Lemma~\ref{irredvarlem} in
Section~\ref{SecCanSpec} will justify that decision.

\section{Reviewing BUILDTREE}\label{SecBuildtree}
A {\em specification of specializations} is a subset $\Sigma$ of
specializations determined by a constructible set of the parameter
space.

BUILDTREE uses {\em reduced specifications} for the segments.
Different definitions have been given in~\cite{Mo02}
and~\cite{MaMo06}. The reduced specification used in release 4 of
the DPGB package described in~\cite{Mo02} does not require $N$ to
be radical nor to obtain a prime decomposition of $N$. In this
approach, when we need to test whether a polynomial in $K[\oa]$
vanishes for $\sigma \in \Sigma$, it is not sufficient to divide
it by $N$. Instead, we must test if it belongs to $\sqrt{\langle N
\rangle}$. But this is simpler than computing the radical and its
prime decomposition. This makes REDSPEC more efficient but does
not ensure all the nice properties that we want to have.
Nevertheless, even if this is a good practical solution, for
theoretical purposes we need to replace the concept of {\em
reduced specification}\footnote{Definition 7 in \cite{MaMo06}.}.

\begin{definition}[Red-specification] \label{redspecdef}
 Given the pair $(N,W)$ of null and not null conditions denote
\[ h= \prod_{w \in W} w \in K[\oa], \ \ \ \hbox{and} \ \ \V(h) =
\bigcup_{w \in W} \V(w) \subset {K'}^m.\]
  They determine a {\em reduced specification of
  specializations {\rm (}red-specification{\rm )}} whenever
  \begin{enumerate}
   \item $N$ is a radical ideal described by its reduced Gr\"obner basis wrt $\succ_{\oa}${\rm )}
   \item $W$ is a set of distinct irreducible polynomials in $K[\oa]$,
   \item Let $N_i$ be the prime components of $\langle N \rangle$ over $K[\oa]$.
   Then $h\not\in N_i$ for all $i$.
   \end{enumerate}
\end{definition}
Note that properties (ii), (iii) of the definition
in~\cite{MaMo06} are simple consequences of Definition
\ref{redspecdef}. Nevertheless property (3) of Definition
\ref{redspecdef} is stronger, and REDSPEC (denoted CANSPEC in
previous papers) is supposed here to verify this new definition of
red-specification.

The segment associated to a red-specification is $S_{(N,W)}=\V(N)
\setminus \V(h)$, and the included specifications are
$\Sigma_{(N,W)}=\{ \sigma_{\oa} : \oa \in \V(N) \setminus \V(h)
\in {{K'}}^m \}$. Let $W=\{\omega_1,\dots,\omega_s\} \subset
K[\oa]$ and $\olam=(\lambda_1,\dots,\lambda_s) \in \Ent_{\geq
0}^s$. Define
\[W(\olam)=\overline{\omega}^{\olam}=\prod_{i=1}^s \omega_i^{\lambda_i}.\]
The set of all non-null polynomials of $K[\oa]$ as a consequence
of $W$ is
\[W^{*}=\{ k W(\olam) \ : \ k\in K,\  \olam \in (\Ent_{\geq 0}^{+})^s \}. \]

\begin{definition}[Reduced polynomial]\label{polred}
A polynomial $f$ is {\em reduced} over the segment $S$ determined
by the red-specification $(N,W)$ if $\overline{f}^N=f$,
$\cont_{\ox}(f)=1$ and $\lc(f) \in W^{*}$.
\end{definition}

\begin{definition}[Good specialization]\label{specwell}
We say that the polynomial $F$ {\em specializes well} to the
reduced polynomial $f$ over the segment $S$ determined by the
red-specification $(N,W)$, if $\overline{F}^N$ and $f$ are
proportional except for non-null normalization, i.e. if
$a\overline{F}^N=bf$ with $a,b \in W^{*}$, (i.e. the coefficients
$a,b$ do not become $0$ on any point of $S$).
\end{definition}

BUILDTREE is a Buchberger-like algorithm. Applied to the ideal $I$
it builds up a rooted binary tree with the following properties:
\begin{enumerate}
  \item  At each vertex $v$ a dichotomic decision is taken about
  the vanishing or not of some polynomial $p(\oa) \in K[\oa]$.
  \item Each vertex is labelled by a list of zeroes and ones; the root label is the empty list.
  At the null child vertex $p(\oa)$ is assumed null and a zero is appended to the
  parent's label, whereas $p(\oa)$ is assumed non-null at the
  non-null son vertex, in which a 1 is appended to the father's label.
  \item  At each vertex $v$, the tree stores $(N_v,W_v)$ and $B_v$,
  where
  \begin{itemize}
   \item [-] $(N_v,W_v)$ determines a reduced specification $\Sigma_v$ of the
   specializations summarizing all the decisions taken in the preceding vertices
   starting from the root.
   \item [-] $B_v$ is reduced wrt $\Sigma_v$ (not faithful) and specializes to a basis of
   $\sigma_{\oa_0}(I)$ for every $\sigma_{\oa_0} \in \Sigma_v$,
   preserving the $\lpp$.
  \end{itemize}
  \item  The set of terminal vertices form a disjoint reduced CGS
  in the sense of Definitions \ref{DefDisCGS} and
  \ref{DefRedCGS}:
  \begin{itemize}
   \item[-] $B_v$ specializes to the
   reduced Gr\"obner basis of $\sigma_{\oa_0}(I)$ for every $\sigma_{\oa_0} \in \Sigma_v$ and has the same
    $\lpp$ set. The polynomials $g$
    in the bases are normalized having $\cont_{\ox}(g)=1$.
   \item [-] The specifications of the set of terminal vertices $t_i$
   determine  subsets $S_{t_i} \subset {{K'}}^m$ forming a partition of
   the whole parameter space ${{K'}}^m$:
   \[\mathcal{X}=\{S_{t_0},S_{t_1},\dots,S_{t_p}\},\]
   and the sets $S_{t_i}$ have characteristic $\lpp$ sets that do not
   depend on the algorithm.
  \end{itemize}
  \item The unique vertex having as label a list of 1
  ($[1,\dots,1]$) corresponds to the {\em generic case} as it is determined
  by only non-null conditions. It does not necessarily contain the whole generic
  case, as we will see next.
\end{enumerate}
Thus the terminal vertices of BUILDTREE form a disjoint reduced
CGS.

\section{Finding a canonical CGS}\label{SecFindCan}
The objective of the paper is to advance in the definition and
computation of a unique (canonical) CGS. In order to reach this
objective we need to obtain an intrinsic family of subsets $S_i$
for the CGS, uniquely determined.

Denote $\Gamma=(C_1,\dots,C_s)$ the disjoint reduced CGS built by
BUILDTREE, i.e. the list of terminal cases $C_i=(B_i,S_i)$. We
shall always set the generic case as the first element of
$\Gamma$. We group them by $\lpp$.
\[\Gamma=((C_{11},\dots,C_{1s_1}),\dots,(C_{k1},\dots,C_{ks_k}))=(\Gamma_1,\dots,\Gamma_k)\]
where the first index denotes the $\lpp$ and, as usual, $C_{11}$
corresponds to the fundamental segment of the generic case.
Obviously the sets in each group $(C_{i1},\dots,C_{is_i})$ are
canonically separated because it cannot exist a common reduced
basis for them, as reduction implies preservation of the $\lpp$.
Thus if it is possible to obtain a unique reduced basis for each
group then we will have a canonical CGS. But even when this is not
possible and one or more groups must be split into several
subgroups forming canonical equivalence classes where each class
admits a common reduced basis, we will have a canonical CGS. Thus
our objective is to obtain this classification.
\begin{conjecture}\label{ConjCan}
Let $\Gamma_i=(C_{i1},\dots,C_{is_i})$ be the set of all segments
of a disjoint reduced CGS having reduced bases with a common
$\lpp_i$. If $C_{ij}$ and $C_{ik}$ admit a common reduced basis
and $C_{ik}$ and $C_{il}$ also, then it exists a common reduced
basis to $C_{ij}$, $C_{ik}$ and $C_{il}$.
\end{conjecture}
If the conjecture is true then we have an equivalent relation
between the segments in $\Gamma_i$ that is independent of the
algorithm and thus shows the existence of the canonical CGS. This
canonical CGS will be also minimal in the sense that it contains
the minimum number of segments of a disjoint reduced CGS. We are
now concerned with the task of giving algorithms to carry out the
task of summarizing the $C_{ij}$ varying $j$ forming the
equivalent class with a unique reduced basis.
\subsection{Using sheaves}\label{SubSecShaves}
Before tackle that task we need to know that in some special cases
we will need to use sheaves. We are indebted to Michael
Wibmer~\cite{Wi06} for the idea of using sheaves for summarizing
some kind of segments, as they are needed for some special
problems. Let us give an example from him.

%\begin{exem}
\noindent {\bf Exemple 1.} Let $I=\langle ax+b,cx+d \rangle$.
Applying BUILDTREE with $\succ_{\ox}\ =\lex(x,y)$ and
$\succ_{\oa}\ =\lex(a,b,c,d)$ we obtain the tree of Figure
\ref{DispgbSheafTree}, that provides the following segments
(ordered by $\lpp$):
\begin{figure}
\begin{center}
\includegraphics{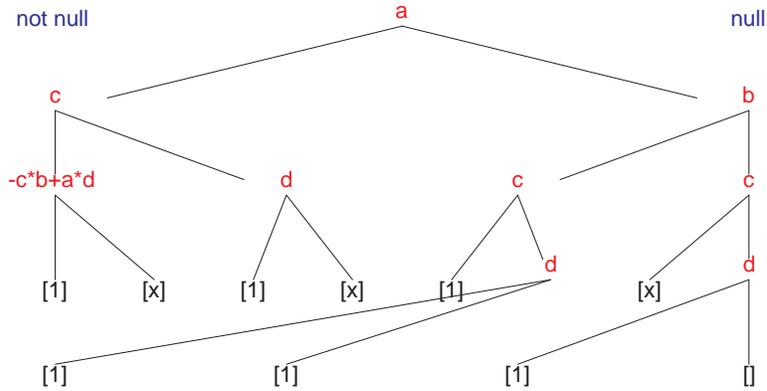}
\caption{\label{DispgbSheafTree} BUILDTREE for $I=\langle
ax+b,cx+d \rangle$.}
\end{center}
\end{figure}

\begin{center}
\begin{tabular}{|l|l|l|l|}
 \hline
 $\lpp$ & basis & null cond. & non-null cond \\
 \hline
 $[1]$ & $[1]$ & $[\ ]$ & $ \{a,ad-cb, c\}$ \\
 $[1]$ & $[1]$ & $[c]$ & $ \{a,d\}$ \\
 $[1]$ & $[1]$ & $[a]$ & $ \{b,c\}$ \\
 $[1]$ & $[1]$ & $[c,a]$ & $ \{b,d\}$ \\
 $[1]$ & $[1]$ & $[d,c,a]$ & $ \{b\}$ \\
 $[1]$ & $[1]$ & $[c,b,a]$ & $ \{d\}$ \\
\hline
 $[x]$ & $[cx+d]$ & $[ad-cb]$ & $\{a,c\}$ \\
 $[x]$ & $[ax+b]$ & $[d,c]$ & $\{a\}$ \\
 $[x]$ & $[cx+d]$ & $[b,a]$ & $\{c\}$ \\
\hline
 $[\ ]$ & $[\ ]$ & $[d,c,b,a]$ & $\{\ \}$ \\
\hline
\end{tabular}
\end{center}

Obviously, the six cases with basis $[1]$ can be summarized into a
single case. We must add the six corresponding segments, and thus
we will need a method to do this in a canonical form. But in any
case, the union of the six segments is intrinsic to the problem
and corresponds to the total generic case having basis $[1]$. You
can see in Section~\ref{SecFurther} how these segments are grouped
in the canonical tree build by the MCCGS algorithm.

The three cases with $\lpp=[x]$ can also be summarized into a
unique basis but now instead of a single polynomial we must use a
sheaf with two polynomials. Effectively, the polynomial $cx+d$
specializes well in the first and third segments with $\lpp=[x]$
and specializes to $0$ in the second segment. And the polynomial
$ax+b$ that forms the reduced basis of the second segment is
proportional (and thus equivalent) to $cx+d$ in the first segment,
but specializes to $0$ in the third segment. The common basis for
the three segments in this case is given by one sheaf
$[\{cx+d,ax+b\}]$ instead by a polynomial: at least one of the two
polynomials in the sheaf specializes well in the union of the
segments whether the other either specializes also or goes to
zero. Thus we see that for our objective we must admit sheaves
also for the bases instead of single polynomials. The three
segments are grouped in the canonical tree.

When a reduced basis of a segment $S$ contains a sheaf, then we
need that, for all $\oa \in S$, at least one of the polynomials in
the sheaf specializes to the corresponding polynomial of reduced
Gr\"obner basis of the specialized ideal and the others either
specialize also, either to it or to $0$.

Thus the canonical CGS for this example will contain only three
segments, namely

\begin{center}
\begin{tabular}{|l|l|l|}
 \hline
 $\lpp$ & basis & sets of pairs $(N,W)$ \\
 \hline
 $[1]$ & $[1]$ & $([\ ],\{a,ad-cb, c\})$, $([c],\{a,d\})$, $([a],\{b,c\})$,\\
       &       & $([c,a],\{b,d\})$,  $([d,c,a],\{b\})$, $([c,b,a],\{d\})$  \\
\hline
 $[x]$ & $[\{cx+d,ax+b\}]$ & $([ad-cb],\{a,c\})$, $([d,c],\{a\})$, $([b,a],\{c\})$ \\
\hline
 $[\ ]$ & $[\ ]$ & $([d,c,b,a],\{\ \})$ \\
\hline
\end{tabular}
\end{center}
%\end{exem}

Sheaves will appear only in over-determined systems with generic
basis $[1]$. For these kind of systems, and for a combination of
the parameter values making compatible the redundance with some
degree of freedom as is the case in the previous example, sheaves
may appear. Nevertheless this is not so for other kind of systems.
An example having a larger sheaf for the basis of one of his
segments is $I=\langle ax+b,cx+d,ex+f\rangle$.

\subsection{Obtaining common reduced bases}
In most common situations where the BUILDTREE CGS presents
multiple segments with the same $\lpp$ it will exist one
subsegment (the most generic one) whose basis already specializes
well in the other segments, and then we only have to pack them.

There are also problems where it does not exist a common basis for
segments having the same $\lpp$. A new example from
Wibmer~\cite{Wi06} shows this situation.

%\begin{exem}
\noindent {\bf Exemple 2.} Consider the following simple system:
$I= \langle u(ux+1), (ux+1)x \rangle$. REBUILDTREE gives the
following GCS with two unique segments having the same $\lpp$.

\begin{center}
\begin{tabular}{|l|l|l|l|}
 \hline
 $\lpp$ & basis & null cond. & non-null cond \\
 \hline
 $[x]$ & $[ux+1]$ & $[\ ]$ & $\{u\}$ \\
 $[x]$ & $[x]$ & $[u]$ & $\{\ \}$ \\
\hline
\end{tabular}
\end{center}

It is easy to convince oneself that it does not exist a common
reduced basis for both segments as the leading term of the generic
segment specializes to $0$ for $u=0$ whether the independent term
is always different from zero.
%\end{exem}

A fourth possibility arises when we have two segments with the
same $\lpp$ sets, characterized by $(B_1,N_1,W_1)$ and
$(B_2,N_2,W_2)$ that do not directly specialize one to the other
by reducing the basis, but nevertheless it can exist a more
generic reduced basis specializing to both. Let us explore that
possibility. We want to test if it exists a basis $B_{12}$ such
that
\[
\begin{array}{rl}
  \hbox{\rm (i)}  & \lpp(B_{12})=\lpp(B_1)=\lpp(B_2).\\
  \hbox{\rm (ii)} & \sigma_{(N_1,W_1)}(B_{12})=B_1$ \hbox{ and } $\sigma_{(N_2,W_2)}(B_{12})=B_2
\end{array}
\]
%\begin{table}[ht]
\fbox{
\begin{minipage}[ht]{5.8in}
 $F \leftarrow \hbox{{\bf GENIMAGE}}(f_1,N_1,W_1,f_2,N_2,W_2)$ \newline
 {\tt Input}: \newline
 \hs\hs $(f_1=\sum a_{\alpha} x^{\alpha},N_1,W_1)$: basis and red-spec of a terminal case \newline
 \hs\hs $(f_2=\sum b_{\alpha} x^{\alpha},N_2,W_2)$: basis and red-spec of a terminal case  \newline
 \hs\hs $L,M\in \Ent$ bounds for the tests. \newline
 {\tt Output}: \newline
 \hs\hs  $F$: when it exists, $F$ returns a polynomial such that $\sigma_{(N_1,W_1)}(F)=f_1$\newline
 \hs\hs\hs\hs  and $\sigma_{(N_2,W_2)}(F)=f_1$ else it returns $F=$ {\bf false}.\newline
 {\bf begin} \newline
 \hs test:= {\bf false}; $F:=$ {\bf false}\newline
 \hs $N:=\hbox{\rm\bf GBEX}(N_1+N_2)$ (returns the Gr\"obner basis and also the matrix $M$\newline
 \hs\hs  expressing the polynomials in $N$ in terms of the polynomials in $N_1$ and $N_2$)\newline
 \hs {\bf for all} $\olam \in \Ent_{\geq 0}^s$,\ $| \olam | \leq L$
 {\bf while} {\bf not} test {\bf do} $w_1=W_1(\olam)$   \newline
 \hs \hs {\bf for all} $\overline{\mu} \in \Ent_{\geq 0}^r$,\ $| \overline{\mu} | \leq M$
 {\bf while} test {\bf do} $w_2=W_2(\overline{\mu})$ \newline
 \hs\hs\hs \{ $HT$ is the index of the leading term) \newline
 \hs\hs\hs $h:=k_1 w_1 a_{HT}-k_2w_2b_{HT}$ \newline
 \hs \hs \hs test := {\bf true} \newline
 \hs\hs\hs {\bf if} $\overline{h}^{N}$ has a factor $Ak_1+Bk_2$
 with $A,B \in K$ {\bf then} \newline
 \hs\hs\hs\hs $k'_1:=B$,\ $k'_2:=-A$ \newline
 \hs\hs\hs\hs {\bf for all} terms $\alpha$ of $f_1$ or $f_2$ {\bf while} test {\bf do}\newline
 \hs\hs\hs\hs\hs $h:=k'_1w_1a_{\alpha}-k'_2w_2b_{\alpha}$ \newline
 \hs\hs\hs\hs\hs $r:=\overline{h}^N$ \newline
 \hs\hs\hs\hs\hs {\bf if} $r \neq 0$ {\bf then} test := {\bf false end if}\newline
 \hs\hs\hs\hs {\bf end for}\newline
 \hs\hs\hs\hs {\bf if} test {\bf then} \newline
 \hs\hs\hs\hs\hs F:=0\newline
 \hs\hs\hs\hs\hs {\bf for all} indices $\alpha$ of terms of $f_1$ or $f_2$  {\bf do}\newline
 \hs\hs\hs\hs\hs\hs $h:=k'_1w_1a_{\alpha}-k'_2w_2b_{\alpha}$ \newline
 \hs\hs\hs\hs\hs\hs $q_i:=\hbox{list of quotients of the exact division } \overline{h}^N$ \newline
  \hs\hs\hs\hs\hs\hs $F:=F+\left(k'_1w_1a_{\alpha} - \displaystyle\sum_{i=1}^{|N|} \sum_{j=1}^{|N_1|}N_{1j} q_i M_{ij}\right)x^{\alpha}$\newline
 \hs\hs\hs\hs\hs {\bf end for}\newline
 \hs\hs\hs {\bf else} \newline
 \hs\hs\hs\hs test := {\bf false} \newline
 \hs\hs\hs {\bf end if}\newline
 \hs\hs {\bf end do} \newline
 \hs {\bf end do} \newline
 {\bf end} \newline
\end{minipage}
}%end fbox
\newline
%\caption{\label{AlgGenimage}}
%\end{table}
in order that both cases can be summarized into a single one
conserving the $\lpp$. Denote $f_1\in B_1$ and $f_2\in B_2$ two
corresponding polynomials with the same $\lpp$. Then we must test
if it exists a $F_{12}$ such that $\sigma_{(N_1,W_1)}(F_{12})=f_1$
and $\sigma_{(N_2,W_2)}(F_{12})=f_2$. For this it must exist
$w_1\in W^{*}_1$, $n_1\in \langle N_1 \rangle \cdot K[\ox]$ and
$w_2\in W^{*}_2$, $n_2\in \langle N_2 \rangle \cdot K[\ox]$ such
that
\[F_{12}=w_1f_1+n_1 = w_2f_2+n_2.\]
Let $f_1=\Sigma_{\alpha} a_{\alpha} x^{\alpha}$,
$f_2=\Sigma_{\alpha} b_{\alpha} x^{\alpha}$ and
$F_{12}=\Sigma_{\alpha} c_{\alpha} x^{\alpha}$. For every index
$\alpha$ of a term in $f_1$ or $f_2$ the coefficients must verify
\[c_{\alpha}=w_1a_{\alpha}+n_{1\alpha} = w_2b_{\alpha}+n_{2\alpha}.\]
with fixed $w_1\in  W^{*}_1$ and  $w_2\in  W^{*}_2$ and
appropriate $n_{1\alpha} \in \langle N_1 \rangle$ and $n_{2\alpha}
\in \langle N_2 \rangle$. This implies
\[w_1a_{\alpha}-w_2b_{\alpha} =n_{2\alpha}-n_{1\alpha} \in \langle
N_1\rangle+\langle N_2 \rangle=\langle N \rangle
\]
that can be solved by testing for all the possible $w_1\in
W^{*}_1$ and $w_2\in W^{*}_2$ when it exists a left hand side
belonging to $\langle N  \rangle$. The semi-algorithm GENIMAGE %of Table \ref{AlgGenimage}
does it and will obtain an $F_{12}$ if it
exists. It is a semi-algorithm in the sense that the possible
choices of $w_1$ and $w_2$ are in fact not finite and the
algorithm must set a bound on the possible total degree ($| \olam
| = \sum_i^s \lambda_i \leq L$) of the terms tested for which no
bound is known. Even if this depends on a luck, and little
combinatorics is used, in practice it does not cause big problems
because this does not occur often and when it does the result is
easily found.

The semi-algorithm is self understanding. Also, when two or more
segments given by red-specifications $(N_1,W_1),\dots, (N_s,W_s)$
have been generalized to a generic basis $B_0$ and we must test if
a new segment $(B,(N,W))$ admits a common pre-image, we can also
use GENIMAGE for each polynomial $f_0\in B_0$ and $f\in B$ taking
for $f_0$ as null and non-null common conditions $(\cap_i^s N_i,
\cap_i^s W_i)$. If GENIMAGE obtains a pre-image it will reduce
well to all the segments. Let us give an example:

%\begin{exem}
\noindent{\bf Exemple 3.} Consider the following example from
Sato-Suzuki~\cite{SaSu03}: $I=\langle
ax^2y+a+3b^2,a(b-c)xy+abx+5c\rangle$ wrt $\succ_{\ox}\ =
\lex(x,y)$,  $\succ_{\oa}\ = \lex(a,b,c)$. REBUILDTREE obtains a
CGS with three segments with basis $[1]$ that obviously can be
directly added, three cases with $\lpp$ set $[y,x]$ that do not
specialize one to the other, and five other segments with distinct
$\lpp$ namely $[y^2, x], [y, x^2], [yx, x^2], [yx^2], [\ ]$.

The question arises for the three segments with $[y,x]$ as $\lpp$.
Let us detail these segments:

\begin{center}
\begin{tabular}{|l|l|l|l|}
 \hline
 $\lpp$ & basis & null cond. & non-null cond \\
 \hline
 $[y,x]$ & $[y,3b^3x-5c]$ & $[a+3b^2]$ & $\{b-c, c, b\}$ \\
 $[y,x]$ & $[a^2y+25, 5x+a]$ & $[b]$ & $\{c,a\}$ \\
 $[y,x]$ & $[25y+3ac^2+a^2,ax+5]$ & $[b-c]$ & $\{c,a\}$ \\
\hline
\end{tabular}
\end{center}
We can verify that none of the bases reduces to the others.
Applying GENIMAGE first to $(B_1,(N_1,W_1))$ and $(B_2,(N_2,W_2))$
and then to $(B_{12}$,$(N_{12}=[b(a+3b^2)],W_{12}=\{c\}))$ and
$(B_3,(N_3,W_3))$ a common reduced basis is found
\[
\begin{array}{lcl}
B_{123} &=& [(25 b c-25 a^3 b-75 b^3 a^2+25 c a^3+75 a^2 b^2 c) y \\
&&-625 a b+1875 c b^2+625 a c+a^2 b c+3 a b^3 c-1875 b^3,\\
&&(b a^2-15 a b+15 a c-9 b^5+9 b^4 c-45 b^3+45 c b^2) x\\
&&-3 b a^2+3 a^2 c+5 a b+27 b^5-27 b^4 c+15 b^3-15 b c^2]
\end{array}
\]
that reduces to the three bases in the respective segments.

%\end{exem}

We have explored four possible situations for a pair (or a
collection) of segments $(B_1,N_1,W_1)$ and $(B_2,N_2,W_2)$ with
the same $\lpp$, namely
\begin{enumerate}
  \item the polynomials of $B_1$ reduce to the polynomials of
  $B_2$ on $(N_2,W_2)$, (most frequent case);
  \item first case does not happen but it exists a pre-image basis
  $B_{12}$ that reduces to both and can be computed by GENIMAGE;
  \item both bases can be summarized using sheaves;
  \item a reduced common basis does not exist.
\end{enumerate}
Table~\ref{AlgDecide} shows the algorithm DECIDE that decides if
two corresponding polynomials of $B_1$ and $B_2$ have a common
pre-image or a sheaf or it does not exist. If
$\overline{S(f_1,f_2)}^{N_2}\neq 0$ then it calls GENIMAGE that
will decide if a pre-image exists or not, but in this case the
result cannot be a sheaf. If $\overline{S(f_1,f_2)}^{N_2}= 0$ as
$\lc(f_1)f_2-\lc(f_2)f_1$ specializes to $0$ in the subset $S_2$,
then $f_1$ specializes either to $f_2$ or to $0$ in $S_2$. Then,
if $\overline{\lc(f_1)}^{N_2} \in W_2^{*}$ then it is always
non-null in $S_2$ and so $f_1$ specializes to $f_2$ and is the
generic polynomial $F$. Else we carry out the symmetric
comparisons and conclusions, and only when the $S$-polynomial
specializes to 0 both in $S_1$ and in $S_2$ and none of the
leading coefficients remain non-null in the other segment we will
have a sheaf $\{f_1,f_2\}$.

\begin{table}[ht]
\fbox{
\begin{minipage}[t]{5.8in}
 $F \leftarrow \hbox{{\bf DECIDE}}(f_1,N_1,W_1,f_2,N_2,W_2)$ \newline
 {\tt Input}: \newline
 \hs\hs $(f_1,N_1,W_1)$: basis and red-spec of a terminal case \newline
 \hs\hs $(f_2,N_2,W_2)$: basis and red-spec of a terminal case  \newline
 {\tt Output}: \newline
 \hs\hs $F$: if it exists a pre-image or a sheaf then $F$ is a polynomial (or sheaf) such that  \newline
 \hs\hs\hs $\sigma_{(N_1,W_1)}(F)=f_1$ and $\sigma_{(N_2,W_2)}(F)=f_2$ else it returns {\bf false}\newline
 {\bf begin} \newline
 \hs {\bf if} $ \overline{S(f_1,f_2)}^{N_2} \neq 0 $ {\bf then} \newline
 \hs\hs $F:=$ {\bf GENIMAGE}$(f_1,N_1,W_1,f_2,N_2,W_2)$\newline
 \hs {\bf else}\newline
 \hs\hs {\bf if} $\overline{\lc(f_1)}^{N_2} \in W_2^{*}$ {\bf then} $F:=f_1$\newline
 \hs\hs {\bf else} \newline
 \hs\hs\hs {\bf if} $ \overline{S(f_1,f_2)}^{N_1} \neq 0 $ {\bf then} \newline
 \hs\hs\hs\hs $F:=$ {\bf GENIMAGE}$(f_1,N_1,W_1,f_2,N_2,W_2)$\newline
 \hs\hs\hs {\bf else}\newline
 \hs\hs\hs\hs {\bf if} $\overline{\lc(f_2)}^{N_1} \in W_1^{*}$ {\bf then} $F:=f_2$\newline
 \hs\hs\hs\hs {\bf else} $F:=\{f_1,f_2\}$ \# (sheaf)\newline
 \hs\hs\hs\hs {\bf end if}  \newline
 \hs\hs\hs {\bf end if}\newline
 \hs\hs {\bf end if}\newline
 \hs {\bf end if}\newline
 {\bf end} \newline
\end{minipage}
}%end fbox
\newline
 \caption{\label{AlgDecide}}
\end{table}

\section{Canonical specifications} \label{SecCanSpec}
The following Lemma plays an important role in the obtention of
canonical specifications of diff-specifications.
%\subsection*{Notation remark}
\begin{lemma}\label{irredvarlem}
Let $K$ be a field of characteristic zero and $K'$ an
algebraically closed extension, $P$ and $Q$ ideals in $K[\oa]$,
$P$ prime and $Q \not\subset P$.  Then, on ${K'}^m$
\[\overline{\V(P) \setminus \V(Q)}=\V(P).\]
\end{lemma}
\begin{proof}
Denote $P'=P\cdot K'[\oa]$ and $Q'=Q\cdot K'[\oa]$ the respective
extensions of $P$ and $Q$ in $K'[\oa]$. To prove the lemma we
follow four steps:\\ \vspace{2mm}
\begin{tabular}{rp{4in}}
  (i) & As $P$ is prime and $Q \not\subset P$, we conclude that $P : Q =
  P$. We leave the proof as an exercise. \\
  (ii) & $(P:Q)'=P':Q'$. See~\cite{ZaSa79}, Vol II, p.
  221. \\
  (iii) & As $P$ is prime, $P'$ is radical. See
  ~\cite{ZaSa79},Vol II, p. 226. \\
  (iv) & Since $K'$ is algebraically closed and $P'$
  is radical,\\
  & \ \ \ $ \overline{\V(P') \setminus \V(Q')}=\V(P' : Q')$. \\
  & See~\cite{CoLiSh92}, Theorem 7, p. 192.
\end{tabular}\\ \vspace{2mm}
Combining these four steps, we obtain
\[\overline{\V(P) \setminus \V(Q)}=\overline{\V(P') \setminus \V(Q')} = \V(P' : Q') = \V((P:Q)') =
\V(P')=\V(P).\]
\end{proof}

 Using Definition~\ref{redspecdef} we can now prove the following
\begin{theorem}\label{Zarclothm}
 Let $(N,W)$ determine a red-specification. Then we have %, over ${K'}^m$ it  results
 \[\overline{\V(N) \setminus \left( \bigcup_{w \in W} \V(w)\right)} =
 \overline{\V(N) \setminus \V(h)} = \V(N)\]
\end{theorem}
\begin{proof}
 Decompose $\sqrt{\langle N \rangle}=\bigcap_i N_i$ into primes in $K[\oa]$, so that
 \[ \V(N) \setminus \V(h)= \left(\bigcup_{i} \V(N_i)\right)
  \setminus \V(h)=\bigcup_{i} \left( \V(N_i)\setminus
  \V(h)\right).
 \]
 As $(N,W)$ determines a red-specification,  $h\not\in
 N_i$ for all $i$, and thus, applying Lemma~\ref{irredvarlem} for each $i$ it
 results
 \[\overline{\V(N) \setminus \V(h)} = \bigcup_{i} \overline{ \V(N_i)\setminus
  \V(h)} = \bigcup_{i} \V(N_i) = \V(N).\]
\end{proof}

We have seen that if Conjecture~\ref{ConjCan} is true it exists an
intrinsic canonical partition of the parameter space $K'^m$ and a
reduced basis for each segment and from the BUILDTREE output the
algorithms DECIDE and GENIMAGE will obtain it. It is apparent the
need of giving a canonical representation of the intrinsic
partition because otherwise we cannot verify its uniqueness, for
example if determined by another algorithm. So we focus now in the
canonical description of the union of red-specifications.
\begin{definition}[Diff-specification]
Given two ideals $N \subset M$ whose associated varieties verify
$\V(N) \supset \V(M)$, they define a {\em diff-specification}
$(N,M)$ describing the subset $S=\V(N) \setminus \V(M)$ of $K'^m$.
\end{definition}
In particular a red-specification $(N,W)$ is easily transformed
into a diff-spe\-ci\-fi\-ca\-tion. Take $h=\prod_{w\in W} w$ and
$M=N+ \langle h \rangle$. Obviously $(N,M)$ is a
diff-specification.

We begin giving a canonical representation of the subsets of a
diff-specification, and then we shall discuss how to add subsets
defined by diff-specifications.

\begin{definition}[Can-specification] \label{canspecdef}
A {\em can-specification} of a subset $C$ is a representation
defined by the set of prime ideals $(N_i,M_{ij})$ varying $i,j$
such that
\begin{equation}\label{canspecform} C=\V(N) \setminus
\V(M)=\bigcup_{i} \left(\V(N_i) \setminus \left(\cup_{j}
\V(M_{ij}\right) \right),
\end{equation}
where ${\mathcal N}=\cap_i N_i$ and ${\mathcal M}_i=\cap_j M_{ij}$
are the prime decompositions over $K[\oa]$ of the radical ideals
${\mathcal N}$ and ${\mathcal M}_i$ respectively, where $ N_i
\subsetneq {\mathcal M}_{ij}$.
\end{definition}
We have the following
\begin{table}[ht]
\fbox{
\begin{minipage}[t]{5.8in}
 $S \leftarrow \hbox{{\bf PRIMEDECOMP}}(N)$ \newline
 {\tt Input}: \newline
 \hs\hs $N$: ideal (representing a variety)\newline
 {\tt Output}: \newline
 \hs\hs  $S=(N_1,\dots,N_k)$: the set of irredundant prime ideals wrt $K$ of the decomposition \newline
 \hs\hs\hs\hs\hs\hs\hs\hs\hs\hs\hs\hs  of $\sqrt{N}=\cap_j N_j$ \newline
 \vspace{0.3cm}

 $Y \leftarrow \hbox{{\bf DIFFTOCANSPEC}}(N,M)$ \newline
 {\tt Input}: \newline
 \hs\hs $N$: the null-condition ideal of the diff-specification\newline
 \hs\hs $M$: the non-null condition ideal of the diff-specification $M \supsetneq N$\newline
 {\tt Output}: \newline
 \hs\hs  $Y=\{(N_i,(\{M_{ij}: 1 \le j \le \ell_i\})) : 1 \le i \le k \}$: the set of prime ideals corresponding \newline
 \hs\hs\hs\hs to the canonical decomposition of $\V(N) \setminus \V(M)$ (Theorem~\ref{canspecthm})\newline
 %\hs\hs\hs\hs  (Theorem~\ref{canspecthm}) .\newline
 {\bf begin} \newline
 \hs $Y=\emptyset$ \newline
 \hs $S:=\hbox{\bf PRIMEDECOMP}(N)$ \newline
 \hs {\bf for} $N_j \in S$ {\bf do} \newline
 \hs \hs {\bf if} $N_j \ne \sqrt{M+N_j}$ {\bf then} \newline
 \hs \hs \hs $T_j:=\hbox{\bf PRIMEDECOMP}(M+N_j)$ \newline
 \hs \hs \hs $Y:=Y \cup_j \{(N_j,T_j)\}$ \newline
 \hs \hs {\bf end if} \newline
 \hs {\bf end for} \newline
 {\bf end} \newline
\end{minipage}
}%end fbox
\newline
\caption{\label{difftocanspecalg} }
\end{table}
\begin{theorem}\label{canspecthm}
\begin{enumerate}
 \item []
  \item  Every set $C=\V(N) \setminus \V(M)\subset {K'}^m$ corresponding to a diff-specification
  admits a can-specification,  %of the form of  formula~(\ref{canspecform}),
  and the algorithm
  {\rm DIFFTOCANSPEC} given in Table~\ref{difftocanspecalg} builds it.
  \item  Over ${K'}^m$, a can-specification verifies
  \[\overline{C}=\overline{\bigcup_{i} \left(\V(N_i) \setminus \left(\cup_{j}
\V(M_{ij}\right) \right)} = \bigcup_i \V(N_i)=\V({\mathcal N}).\]
  \item  The can-specification associated to a set $C$ given by a diff-specification is unique.
  \item  All points in $C \cap \V(N_i)$ are in $\V(N_i) \setminus \left(\cup_{j}
\V(M_{ij})\right)$.
\end{enumerate}
\end{theorem}

\begin{proof}

\begin{enumerate}
  \item Let $\sqrt{N}=\bigcap_{i} N_i$ be the prime decomposition of the
radical ideal $\sqrt{N}$ over $K[\oa]$. Then we have
\[C=\V(N) \setminus \V(M)=\left(\bigcup_{i} \V(N_i)\right)
\setminus \V(M) = \bigcup_{i} \left(\V(N_i) \setminus \V(M+N_i)
\right). \] In this decomposition the variety to be subtracted
from $\V(N_i)$ is contained in it.

 It can happen that $\sqrt{M+N_i}=\langle 1 \rangle$, in which
case nothing is to be subtracted from $\V(N_i)$. It can also
happen that $\sqrt{M+N_i}=N_i$, in which case the term $\V(N_i)
\setminus \V(N_i)$ disappears from the union. The above expression
is simplified and for all the remaining terms we have $ N_i
\subsetneq \sqrt{M+N_i}$

 Let now $\sqrt{M+N_i}=\bigcap_{j} M_{ij}$ be the prime
decomposition of $\sqrt{M+N_i}$ over $K[\oa]$. For each $j$ we
have $N_i  \subsetneq M_{ij}$. The decomposition becomes
\begin{equation}\label{canrepform}
\V(N) \setminus \V(M) = \bigcup_{i} \left(\V(N_i) \setminus
\left(\cup_{j} \V(M_{ij}\right) \right),
\end{equation}
where ${\mathcal N}=\bigcap_i N_i$ and ${\mathcal M}_i=\bigcap_j
M_{ij}$ are the prime decompositions over $K[\oa]$ of the radical
ideals ${\mathcal N}$ and ${\mathcal M}_i$ respectively, proving
part (i) of the theorem. (Observe that the algorithm DIFFTOCANSPEC
is nothing else than the description done in this paragraph).

It should be noted that the prime decompositions in the
computations are performed in $K[\oa]$, as it is the computable
field. Thus these decompositions can split over $K'[\oa]$. In the
same sense we cannot ensure that the varieties $\V(N_i)$ nor
$\V(M_{ij})$ are irreducible over $K^m$ nor over ${K'}^m$ as we
cannot use the Nullstellensatz in $K$. Nevertheless the prime
decompositions are canonically well defined over $K[\oa]$.\\
  \item Using Lemma~\ref{irredvarlem} for each term in the
decomposition given by formula~(\ref{canspecform}) of $C$ we have
\[\overline{C}=\bigcup_i \overline{\V(N_i) \setminus (\cup_j
\V(M_{ij})}= \bigcup_i \V(N_i)=\V({\mathcal N})\] over ${K'}^m$,
proving part (ii) of the theorem.
  \item  Suppose that $C$ admits two
  diff-specifications characterized by the pairs of discriminant
  ideals $(N,M)$ and $(R,S)$ respectively.
   If we denote by
${\mathcal R}=\bigcap_k R_k$ and ${\mathcal S}_{\ell} =
\bigcap_{\ell} S_{k\ell}$ the ideals in the decomposition obtained
from the diff-specification with $R$ and $S$ using the method
described in part (i) of this theorem they will verify
$\overline{C}=\V({\mathcal N})=\V({\mathcal S})$ by part (ii). As
${\mathcal N}$ and ${\mathcal S}$ are radical, they are also
radical over $K'[\oa]$ and thus by the Nullstellensatz they are
both equal to $\I(C)$ over $K'[\oa]$. Thus we have $N_i=R_i$ for
each $i$, as the prime decomposition in $K[\oa]$ is unique.

Next we subtract from each $\V(N_i)$ all the points that are not
in $C$ as they are in $\V(M)$. We have already eliminated the
components of $\V(N)$ that are also in $\V(M)$, so that the points
in $\V(N_i)$ that are not in $C$ are the points of the variety
$\V(M+N_i)\subsetneq \V(N_i)$. Then by the Nullstellensatz
${\mathcal M}_i= \sqrt{M+N_i}$ is the variety ideal
$\I(\V({\mathcal M}_i))$. Carrying out the prime decomposition of
${\mathcal M}_i$ we are done with the canonical decomposition.
Thus the decomposition of $C$ given in part (i) of formula
(\ref{canspecform}) is unique.
 \item This is now obvious as we have subtracted from
 each $\V(N_i)$ all the points in $\V(N_i) \cap \V(M)$.
\end{enumerate}
\end{proof}
Note that the can-specification is canonical but the constructible
sets whose union describes $C$ do not have empty intersection.
Nevertheless this does not affect $C$ itself.

\section{Further developments}\label{SecFurther}

If Conjecture~\ref{ConjCan} is true, it exists a minimal canonical
CGS. The algorithms here described, start from the BUILDTREE CGS
and regroup the segments to obtain the intrinsic segments having a
reduced basis. The result for each segment is of the form:
\[C_i=(B_i,S_i)=(B_i,((N_{i1},W_{i1}),\dots,(N_{ij_i},W_{ij_i})))\]
The subsegments defined by $(N_{ik},W_{ik})$ are described by
red-specifications and thus as a difference of varieties
$\V(N_{ik}) \setminus \V(M_{ik})$ where $h=\prod_{w\in W_{ik}} w$
and $M_{ik}=\langle N_{ik}\rangle +\langle h \rangle$ thus
$(N_{ik},M_{ik})$ is its diff-specification. In order to obtain a
canonical description of the intrinsic partition of our disjoint
reduced CGS it is apparent that we need to add diff-specified sets
in a canonical form. This task is done in~\cite{MaMo07a}.

\begin{figure}
\begin{center}
\includegraphics{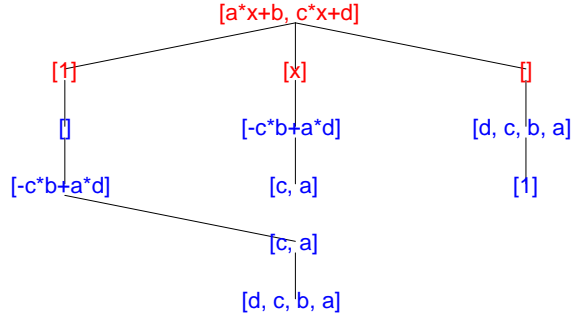}
\caption{\label{FigSheafTree} MCCGS tree for Exemple 1.}
\end{center}
\end{figure}

Let us outline how this works. We cannot assume that the simple
form given by formula~(\ref{canspecform}) will be sufficient. A
more complex constructible set will be formed. There can exist
different canonical forms for describing it, but in any case this
will need prime decomposition of radical ideals. Our canonical
form is given by an even level rooted tree called $P$-tree whose
root defines level $0$. At the nodes there are prime ideals
$P_{i_1,\dots,i_j}$ of $K[\oa]$. These ideals verify
$P_{i_1,\dots,i_j} \subsetneq P_{i_1,\dots,i_j,k}$ for every $k$
and the set of $P_{i_1,\dots,i_j,k}$ for every $k$ are the prime
decomposition of a radical ideal $\mathcal{P}_{i_1,\dots,i_j}$.
The set $C$ defined by the $P$-tree has to be read
\[ C=\bigcup_{i_1} \V(P_{i_1})\setminus
\left(\bigcup_{i_2} \V(P_{i_1i_2}) \setminus \left(\bigcup_{i_3}
\V(P_{i_1i_2i_3}) \setminus \left( \dots \setminus
\bigcup_{i_{2N}} \V(P_{i_1\dots i_{2N}})\right)\right)\right).
\]
Let us advance some results from~\cite{MaMo07a} and illustrate how
is the final output of the MCCGS algorithm for Exemple 1. It is
illustrated by a plot procedure in Figure~\ref{FigSheafTree}, and
the algebraic output summarized in the following table:

\begin{center}
\begin{tabular}{|l|l|l|}
 \hline
 $\lpp$ & basis & segment \\
 \hline
 $[1]$ & $[1]$ & $\Compl^4 \setminus (\V(ad-bc)\setminus (\V(a,c) \setminus \V(a,b,c,d)))$  \\
 \hline
 $[x]$ & $[\{cx+d,ax+b\}]$ & $V(ad-bc) \setminus \V(a,c)$  \\
\hline
 $[\ ]$ & $[\ ]$ & $\V(a,b,c,d)$  \\
\hline
\end{tabular}
\end{center}

It can be seen that the MCCGS algorithm gives a very compact
solution for the problem easy to interpret. This is generally so
for many other problems. It has been successfully applied  to
geometrical theorem discovery~\cite{MoRe07} obtaining very simple
answers for relatively complex problems.

There are two possible lacks coming from the Conjecture and the
semi-algorithm GENIMAGE as we must set artificial bounds to make
it algorithmic. Nevertheless, the use of MCCGS is at least useful
to find examples where the minimal canonical CGS either does not
exist or is not obtained by the actual algorithm, providing
examples to test both the Conjecture and the semi-algorithm.

Finally it must be pointed out that the term order $\succ_{\oa}$
chosen for the computation in $A$ only affect the description by
Gr\"obner bases of the varieties describing the segments but not
to the varieties themselves.

\section*{Acknowledgements} I am indebted to many people. First I
will thank V. Weispfenning for introducing me in this subject, for
his important developments and by his direct help in all the work.
Then I am indebted to R. Sendra for giving me the nice proof of
Lemma~\ref{irredvarlem} that I present here. I am also indebted to
J. Pfeifle for his many helpful comments and frequent discussions.
I am indebted to M. Wibmer for his interesting examples helping me
to understand the sheaf case and other special cases. I would like
to thank P. Viader and Julian Pfeifle for their many helpful
comments and his insightful perusal of the first draft. Finally I
thank my Ph. D. student M. Manubens for helping in all the
discussions and implementation details without which the
theoretical results would not have been possible.


\begin{thebibliography}{99}

% \bibitem[Adams, Boyle(1992)]{AdBo92}
%   W.W.~Adams, A.K.~Boyle.
%   \newblock Some Results on Gr\"obner Bases over Commutative
%   Rings.
%   {\em Jour. Symb. Comp.} (1992) {\bf 13} 473-484.

  \bibitem[Be94]{Be94}
   T.~Becker.
   \newblock On Gr\"obner bases under specialization.
   {\em Appl. Algebra Engrg. Comm. Comput.} (1994) 1--8.

  \bibitem[BeWe93]{BeWe93}
   T.~Becker, V.~Weispfenning.
   \newblock Gr\"obner Bases: A Computational Approach to Commutative
   Algebra.
   \newblock Springer, New-York, (1993).

  \bibitem[Co04]{Co04} M. Coste.
  \newblock Classifying serial manipulators: Computer Algebra and
  geometric insight.
  \newblock Plenary talk. (Personal communication).
  \newblock {\em Proceedings} of EACA-2004 (2004) 323--323.

  \bibitem[CoLiSh92]{CoLiSh92} D.~Cox, J.~Little, D.~O'Shea.
  \newblock Ideals, Varieties and Algorithms.
  \newblock Springer, New-York, (1992). $3^{rd}$ edition (2007).

   \bibitem[De99]{De99} S. Delli\`ere.
   \newblock Triangularisation de syst\`emes constructibles. Application \`a l'\'evaluation dynamique.
   \newblock Th\`ese Doctorale, Universit\'e de Limoges. Limoges, (1995).

   \bibitem[DoSeSt06]{DoSeSt06} A. Dolzmann, A. Seidl, T. Sturm. (2006)
   \newblock REDLOG software in REDUCE http://staff.fim.uni-passau.de/$\sim$\,sturm/

   \bibitem[Du95]{Du95} D. Duval.
   \newblock \'{E}valuation dynamique et cl\^{o}ture alg\'{e}brique en Axiom.
   \newblock {\em Journal of Pure and Applied Algebra} {\bf
   99} (1995) 267--295.

   \bibitem[Em99]{Em99} I. Z. Emiris.
   \newblock Computer Algebra Methods for Studying and Computing Molecular Conformations.
   \newblock {\em Algorithmica}  {\bf 25} (1999) 372-402.

   \bibitem[FoGiTr01]{FoGiTr01}
   E.~Fortuna, P.~Gianni and B.~Trager.
   \newblock Degree reduction under specialization.
   \newblock {\em Jour. Pure and Applied Algebra}, {\bf
   164}:1-2 (2001) 153--164. {\em Proceedings} of MEGA 2000.


   \bibitem[Gi87]{Gi87} P. Gianni.
   \newblock Properties of Gr\"obner bases under specializations.
   \newblock In: EUROCAL'87. Ed. J.H. Davenport, Springer
   LCNS {\bf 378} (1987) 293--297.

   \bibitem[Gom02]{Gom02} T. G\'{o}mez-D\'{\i}az.
   \newblock Dynamic Constructible Closure.
   \newblock {\em Proceedings} of Posso Workshop on Software, Paris, (2000) 73--93.

  \bibitem[GoRe93]{GoRe93} M.J. Gonz\'{a}lez-L\'{o}pez, T. Recio.
  \newblock The ROMIN inverse geometric model and the dynamic
  evaluation method.
  \newblock In: Computer Algebra in Industry.
  \newblock Ed. A.M. Cohen, Wiley \& Sons, (1993) 117--141.

  \bibitem[GoTrZa00]{GoTrZa00} M.J. Gonz\'{a}lez-L\'{o}pez, L. Gonz\'{a}lez-Vega, C. Traverso, A. Zanoni.
  \newblock Gr\"{o}bner Bases Specialization through Hilbert Functions: The Homogeneous Case.
  \newblock  {\em SIGSAM BULL} (Issue 131) {\bf 34}:1 (2000) 1-8.

  \bibitem[GoTrZa05]{GoTrZa05} L. Gonz\'{a}lez-Vega, C. Traverso, A. Zanoni.
  \newblock Hilbert Stratification and Parametric Gr\"{o}bner Bases.
  \newblock {\em Proceedings} of CASC-2005 (2005) 220--235.

   \bibitem[GuOr04]{GuOr04} M. de Guzm\'an, D. Orden.
   \newblock Finding tensegrity structures: geometric and symbolic aproaches.
   \newblock {\em Proceedings} of EACA-2004 (2004) 167--172.

  \bibitem[HeMcKa97]{HeMcKa97} P.~Van~Hentenryck, D.~McAllester and D.~Kapur.
  \newblock Solving polynomial systems using a branch and prune approach.
  \newblock {\em SIAM J. Numer. Anal.} {\bf 34}:2 (1997) 797--827.

   \bibitem[Ka97]{Ka97} M. Kalkbrenner.
   \newblock On the stability of Gr\"obner bases under
   specializations.
   \newblock {\em Jour. Symb. Comp.} {\bf 24}:1 (1997) 51--58.

   \bibitem[Kap95]{Kap95} D.~Kapur.
   \newblock An Approach for Solving Systems of Parametric
   Polynomial Equations.
   \newblock In: Principles and Practices of Constraints Programming.
   \newblock Ed. Saraswat and Van Hentenryck, MIT Press, (1995) 217--244.

   \bibitem[MaMo06]{MaMo06} M.~Manubens, A.~Montes.
   \newblock Improving DISPGB Algorithm Using the Discriminant Ideal.
   \newblock {\em Jour. Symb. Comp.} {\bf 41}:11 (2006) 1245--1263.
   %Short abstract in Proceedings of Algorithmic Algebra and Logic (2005) 159--166.
   %arXiv: math.AC/0601763.

   \bibitem[MaMo07a]{MaMo07a} M.~Manubens, A.~Montes.
   \newblock Minimal Canonical Comprehensive Gr\"obner System.
   \newblock arXiv: math.AC/0611948, (01-12-06). (to be published
   in the {\em Jour. Symb. Comp.}

%   \bibitem[MaMo07b]{MaMo07b}
%   M.~Manubens, A.~Montes.
%   \newblock Rebuilding the BUILDTREE comprehensive Gr\"obner system.
%   \newblock Preprint.

  \bibitem[Mo95]{Mo95}
   A.~Montes.
   \newblock Solving the load flow problem using Gr\"obner bases.
   \newblock {\em SIGSAM Bull.} {\bf 29} (1995) 1--13.

  \bibitem[Mo98]{Mo98}
   A.~Montes.
   \newblock Algebraic solution of the load-flow
     problem for a 4-nodes electrical network.
   \newblock {\em Math. and Comp. in Simul.} {\bf 45} (1998) 163--174.

   \bibitem[Mo02]{Mo02} A.~Montes.
   \newblock New algorithm for discussing Gr\"obner bases with parameters.
   \newblock {\em Jour. Symb. Comp.} {\bf 33}:1-2 (2002) 183--208.

   \bibitem[MoRe07]{MoRe07} A.~Montes, T.~Recio.
   \newblock Automatic discovery of geometry theorems using minimal canonical comprehensive Groebner
   systems.
   \newblock arXiv: math/0703483.

   \bibitem[Mor97]{Mor97} M. Moreno-Maza.
   \newblock Calculs de Pgcd au-dessus des Tours d'\'{E}xtensions
   Simples et R\'{e}solution des Syst\`{e}mes d'\'{E}quations Algebriques.
   \newblock Doctoral Thesis, Universit\'{e} Paris 6, 1997.

   \bibitem[Pe94]{Pe94}
   M.~Pesh.
   \newblock Computing Comprehesive Gr\"obner Bases using MAS.
   \newblock User Manual, Sept. 1994.

    \bibitem[Ry00]{Ry00} M.~Rychlik.
   \newblock Complexity and Applications of Parametric Algorithms of Computational
   Algebraic Geometry.
   \newblock In: Dynamics of Algorithms.
   \newblock Ed. R. del la Llave, L. Petzold, and J. Lorenz.
   \newblock IMA Volumes in Mathematics and its Applications,
   Springer-Verlag,
   {\bf 118} (2000) 1--29.

%   \bibitem[Sa05]{Sa05} Sato, Y., 2005.
%   \newblock Stability of Gr\"obner basis and ACGB.
%   \newblock Proceedings of A3L 2005 (Conference in Honour of the 60th
%   Birthday of V. Weispfenning), eds. A. Dolzmann, A Seidl, T.
%   Sturm. p 223-228. BOD Norderstedt.

   \bibitem[SaSu03]{SaSu03} Sato, Y., Suzuki, A., 2003.
   \newblock An alternative approach to Comprehensive Gr\"obner bases.
   \newblock {\em Jour. Symb. Comp.} {\bf 36}:3-4 (2003), 649-667.

   \bibitem[SuSa06]{SuSa06} Suzuki, A., Sato, Y., 2006.
   \newblock A Simple Algorithm to compute Comprehensive Gr\"obner
   bases. Proceedings of ISSAC 2006, ACM. p 326-331.
   \newblock Implementation in Risa/Asir and Maple
   (http://kurt.scitec.kobe-u.ac.jp/$\sim$sakira/).

%   \bibitem[SaSuNa03]{SaSuNa03}
%   Y. Sato, A. Suzuki, K Nabeshima.
%   \newblock ACGB on Varieties.
%   \newblock {\em Proceedings} of CASC 2003. Passau University, (2003) 313--318.

   \bibitem[Sc91]{Sc91}
   E.~Sch\"onfeld.
   \newblock Parametrische Gr\"obnerbasen im Computeralgebrasystem
   ALDES/SAC-2.
   \newblock Dipl. thesis, Universit\"at Passau, Germany, May 1991.

   \bibitem[Si92]{Si92} W. Sit.
   \newblock An Algorithm for Solving Parametric Linear Systems.
   \newblock {\em Jour. Symb. Comp.} {\bf 13} (1992) 353--394.

  \bibitem[We92]{We92}
   V.~Weispfenning.
   \newblock Comprehensive Gr\"obner Bases.
   \newblock {\em Jour. Symb. Comp.} {\bf 14} (1992) 1--29.

   \bibitem[We03]{We03}
   V.~Weispfenning.
   \newblock Canonical Comprehensive Gr\"obner bases.
   \newblock {\em Proceedings} of ISSAC 2002, ACM-Press, (2002) 270--276.
   {\em Jour. Symb. Comp.} {\bf 36} (2003) 669--683.

   \bibitem[Wi06]{Wi06} M.~Wibmer.
   \newblock Gr\"obner bases for families of affine schemes.
   \newblock arXiv. math/0608019 (2006).

   \bibitem[ZaSa79]{ZaSa79} O. Zariski, P. Samuel. Commutative
   Algebra, 2 volumes. Reprint of the 1958-60 edition, Springer,
   New-York (1979).


%   \bibitem[MaMo05a]{MaMo05a}
%   M.~Manubens, A.~Montes.
%   \newblock Improving DISPGB Algorithm Using the Discriminant
%   Ideal.
%   \newblock Proceedings of Algorithmic Algebra and Logic (2005) 159--166.

%   \bibitem[MaMo05c]{MaMo05c}
%   M.~Manubens, A.~Montes.
%   \newblock Tutorial for the DPGB Library for Discussing Parametric Polynomial
%   Systems.
%   \newblock Submitted to {\em Int. Jour. of Control}.


%  \bibitem[Mo95]{Mo95}
%   A.~Montes.
%   \newblock Solving the Load Flow Problem Using the Gr\"obner Bases.
%   \newblock {\em SIGSAM Bull.}, {\bf 29}:1--13, 1995.

%  \bibitem[Mo98]{Mo98}
%   A.~Montes.
%   \newblock Algebraic solution of the load-flow
%     problem for a 4-nodes electrical network.
%   \newblock {\em Math. and Comp. in Simul.} {\bf 45}:163--174, 1998.

%   \bibitem[Mo02]{Mo02}A.~Montes.
%   \newblock New Algorithm for Discussing Gr\"obner Bases with Parameters.
%   \newblock {\em Jour. Symb. Comp.} {\bf 33}:1-2 (2002) 183--208.

%    \bibitem[Ry00]{Ry00} M.~Rychlik.
%   \newblock Complexity and Applications of Parametric Algorithms of Computational
%   Algebraic Geometry.
%   \newblock In: Dynamics of Algorithms.
%   \newblock Ed. R. del la Llave, L. Petzold, and J. Lorenz.
%   \newblock IMA Volumes in Mathematics and its Applications, Springer-Verlag
%   {\bf 118}:1--29, 2000.

%  \bibitem[We92]{We92}
%   V.~Weispfenning.
%   \newblock Comprehensive Gr\"obner Bases.
%   \newblock {\em Jour. Symb. Comp.} {\bf 14} (1992) 1--29.

%   \bibitem[We02]{We02}
%   V.~Weispfenning.
%   \newblock Canonical Comprehensive Gr\"obner Bases.
%   \newblock {\em Proceedings} of ISSAC 2002. ACM-Press, 270--276, 2002.

%   \bibitem[We03]{We03}
%   V.~Weispfenning.
%   \newblock Canonical Comprehensive Gr\"obner Bases.
%   \newblock {\em Jour. Symb. Comp.} {\bf 36} (2003) 669--683.
%   \newblock Proceedings of ISSAC 2002, ACM-Press, (2002) 270--276.

\end{thebibliography}
\end{document}